\documentclass[a4,12pt,reqno,fleqn]{amsart}
\usepackage{graphicx}
\usepackage{amscd}
\usepackage{xfrac}
\usepackage{mathtools}
\emergencystretch=\maxdimen
\allowdisplaybreaks[1]
\usepackage{array}
\usepackage{xcolor}
\usepackage[all]{xy}


\usepackage{amsmath}
\usepackage{amsthm}
\usepackage{amsfonts,amssymb,hyperref,mathrsfs,cleveref,setspace,enumitem,verbatim}
\usepackage[margin=3.5cm]{geometry}

\DeclareMathAlphabet{\mathpzc}{OT1}{pzc}{m}{it}

\newtheorem{thm}{Theorem}[section]
\newtheorem{cor}[thm]{Corollary}

\newtheorem{propn}[thm]{Proposition}
\theoremstyle{definition}
\newtheorem{defn}[thm]{Definition}
\newtheorem{remark}[thm]{Remark}
\newtheorem{example}[thm]{Example}

\newcommand{\bj}{\perp_{BJ}}
\newcommand{\nbj}{\not\perp_{BJ}}

\newcommand{\ot}{\otimes}

\author[Mohit and R.~Jain]{Mohit and Ranjana Jain}

\address{Mohit, Department of  Mathematics, University of Delhi, Delhi}
\email{mohitdhandamaths@gmail.com}

\address{Ranjana Jain, Department of Mathematics, University of Delhi, Delhi}
\email{rjain@maths.du.ac.in}
\thanks{Research of the first named author is supported by Savitribai Jyotirao Phule Single Girl Child Fellowship vide F.No. 82-7/2022(SA-III)}
\begin{document}
		\title[Local symmetry and smoothness for $C(K,X)$]{Local symmetry and smoothness in the space of vector-valued continuous functions}
	\maketitle
	\textbf{Abstract:} In this article, we characterize the left symmetric points in $C(K,X)$, where $K$ is a compact Hausdorff space  and $X$ is a Banach space. We also provide necessary and sufficient conditions for  the right symmetric points in $C(K,X)$. Further, we identify the smooth points in the space $C_0(K,X)$, $K$ being  locally compact Hausdorff space and $X$ being a Banach space.
		\section{Introduction}
The notion of Birkhoff-James orthogonality (in short B-J orthogonality) is an important tool for the study of  geometry of Banach spaces. An element $x$ in a normed space $X$ over $\mathbb{K}$ is said to be {\it B-J orthogonal} to  $y \in X$ (written as $x\bj y$) if
	$$\|x+\lambda y\|\geq\|x\|,\; \text{for all} \;\lambda\in \mathbb{K}.$$ 
	Unlike the usual orthogonality in the Hilbert spaces, B-J orthogonality fails to be symmetric, that is, $x\bj y$ may not imply $y\bj x$, for $x,y \in X$. However, in order to have a thorough understanding of the geometry of Banach spaces, it is important to discuss the elements which preserve the symmetry, and therefore the concept of left and right symmetric points with respect to the B-J orthogonality was recently introduced.
	In the recent years, many researchers studied the points of local symmetry in various spaces, one can see \cite{sain, ghosh,komuro, saito, vitor,  bose1} for more details. First, let us recall the definition of left (respectively, right) symmetric point.
	\begin{defn}
	An element $x$ in a normed space $X$ is said to be a {\it left symmetric point} (respectively, {\it right symmetric point}) if $x\bj y$ implies $y\bj x$ (respectively, $y\bj x$ implies $x\bj y)$ for all $y\in X$.
	\end{defn}
	
 In 2019, Komuro et. al. \cite[Theorem 3.6, Theorem 4.6]{komuro}  characterized the left and right symmetric points in the space of continuous functions $C(K)$, where $K$ is a compact, Hausdorff space and $X$ is a Banach space. More precisely,  they proved that the left symmetric points are exactly the functions which vanish everywhere except possibly at an isolated point, whereas the right symmetric points are the functions which attain their norms at every point.  We generalize these results for the vector-valued continuous functions.	
	On the other hand, Sundaresan (\cite{sundaresan}) proved that for a compact Hausdorff space $K$ and a Banach space $X$, the uniform norm on $C(K,X)$ is G$\hat{a}$teaux differentiable at a given element $f\in C(K,X)$ if and only if  $M_f =\{k_0\}$ and the norm of $X$ is G$\hat{a}$teaux differentiable at $f(k_0)$, where $M_f:=\{k\in K: \|f(k)\|=\|f\|_{\infty}\}$ is the norm attaining set of $f$. It is well known that the norm of $X$ is G$\hat{a}$teaux differentiable at a non- zero element $x \in X$ if and only if $x$ is smooth, that is, there exists a unique $f\in X^{*}$ satisfying $\|f\|=1$ and $f(x)=\|x\|$. Thus an element $f\in C(K,X)$ is smooth if and only if $M_f=\{k_0\}$ and $f(k_0)$ is a smooth point. We generalize this characterization for a locally compact space $K$ using different techniques.
	
	 In this article, we characterize left symmetric points in the space $C(K,X)$, where $K$ is a compact Hausdorff space, $X$ is a Banach space. In particular, we prove that a non-zero function $f\in C(K,X)$ is left symmetric if and only if $f=f(k_0)\chi_{\{k_{0}\}}$, where $k_{0}$ is an isolated point of $K$ and $f(k_0)$ is a left symmetric point of $X$. This result is also true for a complex Banach space $X$, if, in addition, $K$ is sequentially compact. Further, if $f \in C(K,X)$ is a right symmetric point, then $f$ attains its norm at every point of $K$. The converse is also true if, in addition, $K$ is connected and $f(k)$ is a right symmetric point in $X$ for every $k \in K$. Also, we prove that for a locally compact, Hausdorff space $K$ and a Banach space $X$, an element $f\in C_0(K,X)$ is smooth if and only if $M_f=\{k_0\}$ and $f(k_0)$ is a smooth point of $X$.

	\section{Main Results}
	 We first recall few notations and existing results, and prove some intermediate results which we need for the further development.	 
	For a real normed space $X$ and for $x \in X$, we use  $$ x^{+}:=\{y \in X: \|x+\lambda y\|\geq\|x\|,\;\forall\;\lambda\geq 0 \},$$ $$x^{-} :=\{y \in X:\|x+\lambda y\|\geq\|x\|,\;\forall\;\lambda\leq0 \}.$$
	If $X$ is a complex Banach space, $U=\{u\in \mathbb{C}: |u|=1\;\text{and}\;\text{arg}\ u\in [0,\pi) \}$ and $u\in U$, then denote 
	$$x_{u}^+=\{y\in X: \|x+u\alpha y\|\geq \|x\|,\;\forall\;\alpha\geq0\},$$ 
	$$x_{u}^-=\{y\in X: \|x+u\alpha y\|\geq \|x\|,\;\forall\;\alpha \leq 0\}.$$
	 
	For a real Banach space $X$, Roy et. al. \cite[Theorem 2.1]{roy} characterized the B-J orthogonality of elements in $C(K,X)$ for a compact space $K$ in the following form:
		\begin{thm}\cite[Theorem 2.1]{roy}\label{chrbjcnts}
		Let $K$ be a compact topological space, $X$ be a real normed space and $f,\;g\in C(K,X)$ be non-zero elements. Then $f\bj g$ if and only if there exist $u_1,\;u_2\in M_f$ such that $g(u_1)\in f(u_1)^{+}$ and $g(u_2)\in f(u_2)^{-}$.
	\end{thm}
We first prove a similar result for $C_0(K,X)$, where $X$ is a complex Banach space. The proof is motivated from \cite[Theorem 2.3]{paul1}.
	\begin{thm}\label{dirpaul}
		Let $K$ be a locally compact, sequentially compact, Hausdorff space and $X$ be a complex Banach space. For $f,g\in C_0(K,X)$, $f\bj g$ if and only if for each $u\in U$, there exist $k_u,\;k'_u \in M_f$ such that $g(k_u)\in (f(k_u))_{u}^+$ and $g(k'_u)\in (f(k'_u))_{u}^-$.
	\end{thm}
	\begin{proof}
		Let $f\bj g$ and $u\in U$ be an arbitrary element. 
		Since $K$ is sequentially compact, we can obtain a sequence ${\{k_n\}}$ in $M_{(f+\frac{u}{n}g)}$ converging to, say $k_u \in K$.  
		Then $f(k_n)$ and $g(k_n)$ converge to $f(k_u)$ and $g(k_u)$ respectively, $f$ and $g$ being continuous. Since $f\bj g$, we have
		\begin{align}\label{eq1}
			\bigg\|f+\frac{u}{n}g\bigg\|=\bigg\|f(k_{n})+\frac{u}{n}g(k_{n})\bigg\|\geq \|f\|\geq\|f(k_{n})\|.
		\end{align}
		By taking the limit as $n$ tends to $\infty$ in above expression, we get $k_u\in M_f$. Now consider $\alpha>0$ and let $n_0\in \mathbb{N}$ be such that $\alpha>\frac{1}{{n_0}}$. Then $\|f(k_{n})+u\alpha g(k_{n})\|\geq\|f(k_{n})\|$ for all $n\geq n_0$. 
		If not, then
		\begin{align*}
			\hspace{-0.7cm}	\bigg\|f(k_{n})+\frac{u}{n} g(k_{n})\bigg\| &=\bigg\|\bigg(1-\frac{1}{n\alpha}\bigg)f(k_{n})+\frac{1}{n\alpha}(f(k_{n})+u\alpha g(k_{n}))\bigg\| \\ 
			&\leq \bigg(1-\frac{1}{n\alpha}\bigg)\|f(k_{n})\|+\frac{1}{n\alpha}\|f(k_{n})+u\alpha g(k_{n})\|\\
			& <\bigg(1-\frac{1}{n\alpha}\bigg)\|f(k_{n})\|+\frac{1}{n\alpha}\|f(k_{n})\|\\
			&=\|f(k_{n})\|
		\end{align*} 
		which is a contradiction to \ref{eq1}. Thus, by taking the limits, we get $\|f(k_u)+u\alpha g(k_u)\|\geq\|f(k_u)\|$, so that $g(k_u)\in (f(k_u))_{u}^+$. Similarly, the existence of $k'_u \in M_f$ can be proved by taking a convergent sequence in $M_{(f-\frac{u}{n}g)}$.
		
		For the converse, let $\lambda\in \mathbb{C}$ be arbitrary. Then, we can write $\lambda=u\alpha$ for some $u\in U$ and $\alpha\in \mathbb{R}$. Therefore, if $\alpha\geq0$,
		$$\|f+\lambda g\|\geq \|f(k_u)+u\alpha g(k_u)\|\geq\|f(k_u)\|=\|f\|,$$
		  and if $\alpha<0$, 
		  $$\|f+\lambda g\|\geq \|f(k'_u)+u\alpha g(k'_u)\|\geq\|f(k'_u)\|=\|f\|.$$ 
		  Thus, $f\bj g$ and hence the result.		
	\end{proof}
	Recently, Martin et. al. \cite[Theorem 4.3]{martin} characterized the B-J orthogonality in $C(K,X)$ in terms of directional orthogonality, where $K$ a compact space and  $X$ is a Banach space. Recall that for $x,y\in X$ and $t\in T:=\{t\in \mathbb{K}: |t|=1\}$, $x$ is said to be {\it orthogonal to $y$ in the direction of $t$} (written as $x\perp_{t}y$) if $\|x+\alpha ty\|\geq \|x\|\;\forall\;\alpha\in \mathbb{R}$.		
	\begin{thm}\cite[Theorem 4.3]{martin}\label{chrbjcntsc}
		Let $K$ be a compact Hausdorff space and $X$ be a Banach space. Let $f,g\in C(K,X)$ be such that $M_f$ is connected. Then $f\bj g$ if and only if for each $t\in T$, there exists $k_t\in M_f$ such that $f(k_t)\perp_{t}g(k_t)$. In particular, when $X$ is real, $f\bj g$ if and only if there exists $k_0\in M_f$ such that $f(k_0)\bj g(k_0)$.
	\end{thm}	
 We generalize this result for a locally compact space $K$. The proof goes on the similar lines with slight modifications.
	\begin{propn}\label{lcchrcnts}
		Let $K$ be a locally compact Hausdorff space and $X$ be a Banach space. Let $f,\;g\in C_0(K,X)$ be such that $M_f$ is connected. Then $f\bj g$ if and only if for each $t\in T$ there exists $k_t\in M_f$ such that $f(k_t)\perp_{t}g(k_t)$.
	\end{propn}
	\begin{proof}
		Let $f\bj g$ and $t\in T$ be an arbitrary element. Consider the set $$A=\{\phi(g): \phi\in \text{Ext}(B_{C_0(K,X)^{*}}),\; \phi(f)=\|f\|\}.$$  By \cite[Chapter II, Theorem 1.1]{singer}, $0\in \text{conv}(A)$ as $f\bj g$. From \cite[Lemma 3.3]{bruno} and \cite[Corollary 3]{wolfgang}, we know that $$\text{Ext}(B_{C_0(K,X)^{*}})=\{x^{*}\ot \delta_{k}:x^{*}\in \text{Ext}(B_{X^*}),\;k\in K\},$$ where $x^{*}\ot \delta_{k}:C_0(K,X)\rightarrow \mathbb{K}$ is defined as $(x^{*}\ot \delta_{k})(f)=x^*(f(k))$. Observe that, if $x^*\in B_{X^*}$ satisfies $x^*(f(k))=\|f\|$ for $k\in K$, then $k\in M_f$. Therefore, 
		$$A=\{x^*(g(k)): k\in M_f,\; x^*\in \text{Ext}(B_{X^{*}}),\; x^*(f(k))=\|f\|\}.$$
		 Now, consider the set $$B=\{x^*(g(k)):k\in M_f,\; x^*\in S_{X^{*}},\; x^*(f(k))=\|f\|\}.$$ Clearly, $0\in \text{conv}(B)$ as $A\subseteq B$. Also, by \cite[Lemma 4.4]{martin} we have that, $B$ is a connected set and hence by \cite[Lemma 2.7]{martin}, for given $t\in T$, there exist $k_t\in M_f$ and $x^*\in S_{X^*}$ such that $x^*(f(k_t))=\|f\|=\|f(k_t)\|$ and \text{Re}$(tx^{*}(g(k_t)))=0$. Thus, if $X$ is a real Banach space, $f(k_t)\perp_{t}g(k_t)$ \cite[Theorem 2.1.15]{mal}. If $X$ is a complex Banach space, then by \cite[Theorem 4]{bagchi}, $tf(k_t)\perp_{t} tg(k_t)$ which further implies $f(k_t)\perp_{t}g(k_t)$. 
		
		For the converse, let $\lambda\in \mathbb{K}$ be arbitrary. Write $\lambda=\alpha t$ where $\alpha\in \mathbb{R}$ and $t\in T$. Then $$\|f+\lambda g\|=\|f+\alpha tg\|\geq \|f(k_t)+\alpha tg(k_t)\|\geq \|f(k_t)\|=\|f\|.$$ Thus, $f\bj g$ and hence the result.
	\end{proof}
	
	With all the ingredients prepared we are now ready to prove the main results. We first characterize the left symmetric points in $C(K,X)$.
	
	\begin{thm}\label{leftcntsreal}
		Let $K$ be a compact Hausdorff space and $X$ be a real Banach space. A non-zero element $f\in C(K,X)$ is a left symmetric point if and only if $f=f(k_0)\chi_{\{k_{0}\}}$, where $k_{0}$ is an isolated point of $K$ and $f(k_0)$ is a non-zero left symmetric point of $X$. If $X$ is a complex Banach space then the characterization holds if, in addition, $K$ is sequentially compact.
	\end{thm}
	\begin{proof} Suppose $f\in C(K,X)$ is of the form $f=f(k_0)\chi_{\{k_{0}\}}$ for some isolated point $k_0\in K$, where $f(k_0)$ is a non-zero left symmetric in $X$. Let $g\in C(K,X)$ be a non-zero element such that $f\bj g$. Note that $M_f=\{k_0\}$.
		
		Case (1): $X$ is real Banach space:
		By \Cref{chrbjcnts}, we have $f(k_0)\bj g(k_0)$. Since $f(k_0)$ is a left symmetric point therefore $g(k_0)\bj f(k_0)$.	If $k_0\in M_g$, then $g\bj f$ by \Cref{chrbjcnts}. If $k_0\not\in M_g$, then since $f(k)=0$ for all $k\neq k_0$ and $M_g\neq\emptyset$, there exist $k_1,\;k_2\in M_g$ such that  $f(k_1)\in g(k_1)^{+}$ and $f(k_2)\in g(k_2)^{-}$ so that $g\bj f$. Hence, $f$ is left symmetric.
		
		Case (2): $X$ is complex Banach space: By \Cref{chrbjcntsc}, for each $t\in T$, $f(k_0)\perp_{t}g(k_0)$ and consequently $f(k_0)\bj g(k_0)$. This gives $g(k_0)\bj f(k_0)$ as $f(k_0)$ is a left symmetric point. Now, if $k_0\in M_g$ then $$\|g+\lambda f\|\geq \|g(k_0)+\lambda f(k_0)\|=\|g(k_0)\|=\|g\|$$ and if $k_0\not\in M_g$, then for $k'\in M_g$ and $\lambda\in \mathbb{C}$, we have $$\|g+\lambda f\|\geq \|g(k')+\lambda f(k')\|=\|g(k')\|=\|g\|.$$ Thus, $g\bj f$.

		Conversely, suppose that $f\in C(K,X)$  is a non-zero left symmetric point. Let $k_{0}$ be an element of $M_f$. We first claim that $f(k)=0$ for all $k\neq k_0$. Let, if possible, there exists $k_1\in K\setminus\{k_0\}$ such that $f(k_1)\neq0$. Let $V \subseteq K$ be an open set such that  $k_1 \in V$ and $k_0\not\in V$. By Urysohn's lemma, there exists a continuous function $h:K\rightarrow [0,1]$ such that $h(k_1)=1,\;h(V^c)=0$ and support of $h$ is compact. Define $g:K\rightarrow X$ as $g(k)=h(k)f(k)$, for all $k \in K$. Clearly, $g$ is a non-zero continuous function as $g(k_1)\neq0$. Now, $g(k_0)=0$ implies that $g(k_0)\in\;f(k_0)^{+}$ and $g(k_0)\in\;f(k_0)^{-}$ (when $X$ is real). Also, for each $u\in U$,\; $g(k_0)\in\;f(k_0)_{u}^{+}$ and $g(k_0)\in\;f(k_0)_{u}^{-}$ (when $X$ is complex). Since $k_0\in M_f$, by \Cref{chrbjcnts} and \Cref{dirpaul}, $f\bj g$ so that  $g\bj f$, $f$ being left symmetric. 
		
		If $X$ is real Banach space, then by \Cref{chrbjcnts}, there exist $k'_1,\;k'_2\in M_g$ such that $f(k'_1)\in g(k'_1)^{+}$ and $f(k'_2)\in g(k'_2)^{-}$. Now, $k'_2\in M_g$ implies that $f(k'_2)\neq0$ and $h(k'_2)>0$. Thus, for $\lambda=-h(k'_2)$, we have that $\|g(k'_2) + \lambda f(k'_2)\|=0 < \|g(k'_2)\|$, which is a contradiction. 
		
		If $X$ is complex Banach space, then by \Cref{dirpaul}, there exists $k'\in M_g$ such that $f(k')\in (g(k'))_{1}^{-}$ that is $\|g(k')+\alpha f(k')\|\geq \|g(k')\|\;\forall\;\alpha<0$. But, for $\alpha=-h(k')<0$, we have  $\|g(k')+\alpha f(k')\|=0<\|g(k')\|$, which is a contradiction. 
		
		Hence, in both the cases, $f(k)=f(k_0)\chi_{\{k_0\}}$. 
		
		Finally, we show that $f(k_0)$ is a left symmetric point of $X$. If not, then there exists $0\neq x\in X$ such that $f(k_0)\bj x$ but $x\nbj f(k_0)$. Define $g\in C(K,X)$ as $g=x\chi_{\{k_0\}}$. Clearly $M_{g}=\{k_0\}=M_{f}$ and $f(k_0)\bj g(k_0)$ where as $g(k_0)\nbj f(k_0)$. Thus, $f\bj g$ but $g\nbj f$, by \Cref{chrbjcnts} (when $X$ is real) and by \Cref{chrbjcntsc} (when $X$ is complex). This is a contradiction to the fact that $f$ is a left symmetric point and hence $f(k_0)$ is left symmetric.
	\end{proof}
	\begin{cor}
		Let $K$ be a compact Hausdorff metric space and $X$ be a Banach space.	If either $K$ has no isolated points or $X$ has no non-zero left symmetric points,  then zero is the only left symmetric point in $C(K,X)$.
	\end{cor}
	
	In \cite[Theorem 4.6]{komuro}, Komuro et. al. proved that for a compact Hausdorff space $K$, a function $f\in C(K)$ with $\|f\|_{\infty} =1$ is right symmetric if and only if $f$ is unimodular. Using few of their techniques, we generalize this result for vector-valued functions over the real field which also provide another proof of \cite[Theorem 4.6]{komuro} (for the real-valued functions). 
	\begin{thm}\label{conti}
		Let $K$ be a compact, Hausdorff space and $X$ be a real Banach space. If $f\in C(K,X)$ is a non-zero right symmetric point then $M_f=K$.
	\end{thm}
	\begin{proof}
		Without loss of generality, we may assume $\|f\|=1$. We first claim that $f(k)\neq 0$, for all $k\in K$. Suppose $f(k_0)=0$ for some $k_0\in K$. Let $K_1=\{k\in K: \|f(k)\|\geq\frac{1}{2}\}$. By Urysohn's lemma, there exists a continuous function $h:K\rightarrow[0,1]$ such that $h(k_0)=1$ and $h\arrowvert_{K_1}=0$. Let $x\in X$ be any arbitrary element of norm one and define $g:K\rightarrow X$ as $g(k)=h(k)x+f(k)(1-h(k))$. Clearly, $g\in C(K,X)$ and $\|g(k_0)\|=1$. Observe that $\|g\|=1$ as for an arbitrary $k\in K$, if $k\in K_1$, then $\|g(k)\|=\|f(k)\|\leq1$, and if $k\not\in K_1$, then
		\begin{align*}
			\|g(k)\|&\leq\|x\||h(k)|+\|f(k)\||1-h(k)|\\
			&\leq h(k)+\frac{1}{2}(1-h(k))\\
			&\leq1.	
		\end{align*}
		Now, for any $\lambda\in \mathbb{R}$, 
		$$	\|g+\lambda f\|\geq\|g(k_0)+\lambda f(k_0)\|=\|g(k_0)\|=\|g\|,$$
		so that $g \bj f$. However,  for $k\in K_1$,  $$\|f(k)-\frac{1}{2}g(k)\|=\frac{\|f(k)\|}{2}<1$$
		and for $k\not\in K_1$,
		\begin{align*}
			\|f(k)-\frac{1}{2}g(k)\|&=\big\|f(k)-\frac{1}{2}\big(h(k)x+f(k)(1-h(k))\big)\big\|\\
			&<\frac{1}{2}+\frac{\|h(k)x\|}{2}+\frac{(1-h(k))}{2}\\
			&=1.
		\end{align*}
		Thus, $\|f-\frac{1}{2}g\|<1=\|f\|$ so that $f\nbj g$. This is a contradiction to the fact that $f$ is a right symmetric point. Thus $f$ is non-zero everywhere.
		
		Next, we claim that $M_f=K$. Let, if possible, there exist $k_0\in K$ such that $ k_0 \notin M_f$. By Urysohn's lemma, we have continuous maps $h,\;h':K\rightarrow [0,1]$ such that $h(k_0)=0,\;h\arrowvert_{M_f}=1$ and $h'(k_0)=1,\;h'\arrowvert_{M_f}=0$. Further, since $K$ is a normal space being compact Hausdorff and $M_f$ is a closed subset, therefore there exist two disjoint open subsets $U_1$ and $U_2$ of $K$ such that $\{k_0\}\subseteq U_1$ and $M_f\subseteq U_2$. Also, $M_f$ is compact, so by Urysohn's lemma, there exist continuous functions  $h_1,h_2:K\rightarrow [0,1]$ such that $h_{1}\arrowvert_{M_f}=1$, $h_{1}(U_2^{c})=0$ and $h_{2}(k_0)=1$ and $h_{2}(U_1^{c})=0$. Define $g:K\rightarrow X$ as $g(k)=f(k)h(k)h_{1}(k)-h'(k)h_{2}(k)\frac{f(k_0)}{\|f(k_0)\|}$. It is clear that $g\in C(K,X)$. Observe that $\|g\|=1$. For this, let $k\in K$ be an arbitrary element. Then,
		$$g(k)=
		\begin{cases}
			-h'(k)h_{2}(k)\frac{f(k_0)}{\|f(k_0)\|} & \text{if}\;k\in U_1\subseteq U_2^{c}\\
			f(k)h(k)h_{1}(k) & \text{if}\;k\in U_2\subseteq U_1^{c}\\
			0 & \text{if}\;k\not\in U_{1}\cup U_2.
		\end{cases}
		$$
		Thus, we have $\|g(k)\|\leq1$ and $\|g(k_0)\|=1$, which gives $\|g\|=1$. Note that $M_f\subseteq M_g$ as for $k\in M_f,\ g(k)= f(k)h(k)h_{1}(k) = f(k)$.
		Now, for $\lambda\geq0$ and for any $k_1\in M_f$, 
		$$\|g(k_1)+\lambda f(k_1)\|=\|f(k_1)\||1+\lambda|=1+\lambda\geq1=\|g(k_1)\|$$
		and for $\lambda\leq0$,
		$$\|g(k_0)+\lambda f(k_0)\|=|-1+\lambda\|f(k_0)\||=1-\lambda\|f(k_0)\|\geq1=\|g(k_0)\|.$$ 
		Thus, $f(k_1)\in g(k_1)^{+}$ and $f(k_0)\in g(k_0)^{-}$ so that $g\bj f$, by \Cref{chrbjcnts}. However, for any $k\in M_f$, $$\big\|f(k)-\frac{1}{2} g(k)\big\|=\frac{\|f(k)\|}{2} <\|f(k)\|.$$
		Thus, $g(k)\notin f(k)^{-}$ and again by \Cref{chrbjcnts}, $f\nbj g$ which is a contradiction to the fact that $f$ is a  right symmetric point. This completes the proof.
	\end{proof}
	
	The converse of the above result is not true.
	\begin{example}
		Let $K=[0,1]\cup\{2\}$ and $X= (\mathbb{R}^2, \|\cdot\|_{\max})$. Define $f\in C(K,X)$ as $f=(1,1)\chi_{[0,1]}+(1,0)\chi_{\{2\}}$ so that $M_{f}=K$. However, $f$ is not a right symmetric point. To see this, consider $g\in C(K,X)$ as $g=(\frac{1}{2},1)\chi_{K}$. Then, for any $\lambda\in \mathbb{R}$, 
		$$\|g+\lambda f\|\geq\|g(2)+\lambda f(2)\|\geq1=\|g\|$$ and $$\bigg\|f-\frac{1}{2}g\bigg\|=\frac{3}{4}<1=\|f\|.$$ Thus, $g\bj f$ but $f\nbj g$ .  
	\end{example}
	
	However, a partial converse of \Cref{conti} is true as shown in the next result. 
	\begin{thm}\label{final}
		Let $K$ be a compact, connected Hausdorff space and $X$ be a real Banach space. Let $f\in C(K,X)$ be such that $M_f=K$ and $f(k)$ is right symmetric point in $X$ for all $k\in\;K$. Then $f$ is  right symmetric.
	\end{thm}
	\begin{proof}
		Let $g\in C(K,X)$ be such that $g\bj f$. If $g(k)=0$ for some $k\in K$, then 
		$$ \|f + \lambda g\| \geq \|f(k)\| = \|f\|, $$
		since $M_f=K$, so that $f\bj g$. Let $g(k)\neq 0$ for all $k\in K$. Define $g'\in C(K,X)$ as $g'(k)=\frac{g(k)}{\|g(k)\|}$, then $\|g'\|=1$ and $M_{g'}=K$.  Now, $g\bj f$ therefore, by \Cref{chrbjcnts}, there exist $k_1\;k_2\in M_g$ such that $f(k_1)\in g(k_1)^{+}$ and $f(k_2)\in g(k_2)^{-}$. For any $\lambda \geq 0$,  $$\|g'(k_1)+ \lambda f(k_1)\| = \frac{1}{\|g(k_1)\|}\|g(k_1)+ (\lambda \|g(k_1)\|)f(k_1)\|\geq1=\|g'(k_1)\|,$$ 	 
		so that $f(k_1)\in g'(k_1)^{+}$, and similarly $f(k_2)\in g'(k_2)^{-}$ and hence $g'\bj f$. Since $M_{g'}$ is connected, by \Cref{chrbjcntsc} there exists $k_0\in K$ such that $g'(k_0)\bj f(k_0)$ which implies $f(k_0)\bj g'(k_0)$ as $f(k_0)$ is a right symmetric point in $X$. By homogeneity property, $f(k_0)\bj g(k_0)$. Since  $M_f=K$, for any $\lambda \in \mathbb{R}$, we have $$\|f+\lambda g\|\geq \|f(k_0)+\lambda g(k_0)\|\geq \|f(k_0)\|=\|f\|.$$ Hence, $f\bj g$. This proves the claim.	
	\end{proof}
Lastly, we characterize the smooth points in $C_0(K,X)$. 

	
	\begin{thm}
		Let $K$ be a  locally compact Hausdorff space, $X$ be a Banach space. An element $f\in C_0(K,X)$ is smooth if and only if $M_f=\{k_0\}$ for some $k_0 \in K$ and $f(k_0)$ is a smooth point in $X$.
	\end{thm}
	\begin{proof}
	Let $f$ be smooth and $k_0\in M_f$. Let, if possible, there exists $k_1\in K\setminus\{k_0\}$ such that $\|f(k_1)\|=\|f\|$. In light of Hahn Banach theorem, let $F_0, F_1\in X^{*}$ such that $\|F_0\|=1=\|F_1\|$, $F_0(f(k_0))=\|f(k_0)\|$ and $F_1(f(k_1))=\|f(k_1)\|$. For $i=0,1$, define $\phi_i:C_0(K,X)\rightarrow \mathbb{K}$ as $\phi_i(g)=F_i(g(k_i))$. Clearly, $\phi_i\in C_0(K,X)^{*}$ and $\phi_i(f)=\|f\|$  as $k_i\in M_f$. Since $\|F_i\|=1$, therefore $\|\phi_i\|=1$. We first claim that $\phi_0\neq \phi_1$. For this, let $V$ be an open set such that $k_0\in V$ and $k_1\not\in V$. Define $p:K\rightarrow X$ as $p(k)=f(k)h(k)$, where $h:K\rightarrow [0,1]$ is a continuous function with compact support for which $h(k_0)=1$ and $h(V^c)=0$. Clearly, $p\in C_0(K,X)$ with $\phi_0(p)=\|f(k_0)\|=\|f\|$ and $\phi_1(p)=0$. Thus, $\phi_0$ and $\phi_1$ are two different support functionals corresponding to $f$ which contradicts the hypothesis that $f$ is a smooth point and this proves the claim.
	
	 Next, we claim that $f(k_0)$ is a smooth point. If not, then there exist $G_1, G_2\in X^{*}$ such that $\|G_1\|=1=\|G_2\|$ and $G_1(f(k_0))=\|f(k_0)\|=G_2(f(k_0))$. For $i=1,2$, define $\psi_i \in (C_0(K,X))^*$ as $\psi_i(g)=G_i(g(k_0))$ for all $g\in C_0(K,X)$. Since $G_1$ and $G_2$ are support functionals corresponding to $f(k_0)$, one can easily check that $\psi_1, \;\psi_2$ both are support functionals corresponding to $f$. Also, $G_1\neq G_2$ implies that there exists $y\in X$ such that $G_1(y)\neq G_2(y)$. Therefore, for $p'\in C_0(K, X)$ defined as $p'(k)=yh(k)$, we have $\psi_1(p')\neq \psi_2(p')$ which contradicts the smoothness of $f$. This completes the proof.
		
	Conversely, suppose that  $M_f=\{k_0\}$ and $f(k_0)$ is a smooth point in $X$. In view of \cite[Theorem 2.3.2, Remark 2.3.4]{mal}, it is sufficient to prove that right additivity holds at $f$. For this, let $g,h\in C_0(K,X)$ such that $f\bj g$ and $f\bj h$. Since $M_{f}=\{k_0\}$, therefore by \Cref{lcchrcnts}, for each $t\in T$ we have $f(k_0)\perp_{t}g(k_0)$ and $f(k_0)\perp_{t}h(k_0)$. This gives $f(k_0)\bj g(k_0)$ and $f(k_0)\bj h(k_0)$ which further implies that $f(k_0)\bj (g(k_0)+h(k_0))$ as $f(k_0)$ is a smooth point. Therefore, $$\|f+\lambda (g+h)\|\geq \|f(k_0)+\lambda (g+h)(k_0)\|\geq \|f(k_0)\|=\|f\|.$$ Thus, $f\bj (g+h)$.		
	\end{proof}

	\begin{remark}
		It would be interesting to know whether $f(k)$ is right symmetric for every $k \in K$ if $f$ is right symmetric in $C(K,X)$.  For some specific cases it happens to be true. One instance can be seen in the space $C(K,c_{00})$, where  $c_{00}$ is the space of eventually zero real sequences.  By \cite[Theorem 2.14]{bose1}, it is well known that $c_{00}$ does not possess any non-zero right symmetric point. In fact, the same is true for the space $C(K,c_{00})$.  For this, consider a non-zero element $f\in C(K,c_{00})$ of norm one.  Let us write $f(k)=(a_{1},a_{2},...a_{n(k)},0,0,... )$.  Define $g:K\rightarrow c_{00}$ as $g(k)=f(k)+e_{n(k)+1}$, then clearly $g$ is a continuous function of norm one. For any $k\in K$ and $\lambda\in\mathbb{R}$, observe that $\|g(k)+\lambda f(k)\|\geq1$ and hence $\|g+\lambda f\|\geq \|g\|$, but $\|f(k)-\frac{1}{2} g(k)\|=\max\{\frac{|a_1|}{2},...,\frac{|a_{n(k)}|}{2},\frac{1}{2}\}\leq\frac{1}{2}$ and hence $\|f-\frac{1}{2}g\|<\|f\|$. Thus, $g\bj f$ and $f\nbj g$, so that $f$ cannot be a right symmetric point.	
	\end{remark}

\end{document}